\documentclass[11pt,twoside]{article}
\usepackage{amssymb,mathrsfs,cite,color}
\usepackage{amsfonts,amsmath,latexsym}

\usepackage{indentfirst}
\usepackage[amsmath,thmmarks]{ntheorem}
\numberwithin{equation}{section}
\theoremheaderfont{\normalfont\rmfamily\bf}
\theorembodyfont{\normalfont\it }
\newtheorem{theorem}{\indent Theorem}[section]
\newtheorem{lemma}{\indent Lemma} [section]

\newtheorem{definition}{\indent Definition} [section]

\newtheorem{thm}{\indent Theorem}
 
  \theoremstyle{nonumberplain}

\theorembodyfont{\normalfont \rm }
 \theoremsymbol{$\square$}
  \newtheorem{proof}{\indent Proof}

\allowdisplaybreaks

\def\rr{\mathbb{R}}
\def\rn{{\rr}^n}

\begin{document}

\title{\bf Characterization of Lipschitz spaces via
commutators of the Hardy-Littlewood maximal function \footnote{
Supported by the National Natural Science Foundation of
China (Grant Nos. 11571160 and 11471176).}
  }
\author{Pu Zhang   
         \\
 {\small\it Department of Mathematics, Mudanjiang Normal
University, Mudanjiang 157011, P. R. China}\\
 {\small\it  E-mail: puzhang@sohu.com 
               }
    }
\date{ }
\maketitle

\begin{center}\begin{minipage}{14.5cm}

{\bf Abstract}~  Let $M$ be the Hardy-Littlewood maximal function
and $b$ be a locally integrable function. Denote by $M_b$ and
$[b,M]$ the maximal commutator and the (nonlinear) commutator of
$M$ with $b$. In this paper, the author consider the boundedness
of $M_b$ and $[b,M]$ on Lebesgue spaces and Morrey spaces when $b$
belongs to the Lipschitz space, by which
some new characterizations of the Lipschitz spaces are given.

{\bf Keywords}~  Hardy-Littlewood maximal function; commutator;
Lipschitz space; Morrey space.

{\bf MR(2010) Subject Classification}~ 42B25, 42B20, 42B35, 46E30

\end{minipage}
\end{center}

\medskip

\section{Introduction and Results}

Let $T$ be the  classical singular integral operator, the commutator
$[b,T]$ generated by $T$ and a suitable function $b$ is given by
\begin{equation} \label{equ.1.1}    
[b,T]f = bT(f)-T(bf).
\end{equation}

A well known result due to Coifman, Rochberg and Weiss \cite{crw}
(see also \cite{j}) states that $b\in {BMO(\rn)}$ if and only if the
commutator $[b,T]$ is bounded on $L^p(\rn)$ for $1<p<\infty$. In
1978, Janson \cite{j} gave some characterizations of the Lipschitz
space ${\dot{\Lambda}}_{\beta}(\rn)$ (see Definition \ref{lip}
below) via commutator $[b,T]$ and proved that $b\in
{\dot{\Lambda}_{\beta}(\rn)} (0<\beta<1)$ if and only if $[b,T]$ is
bounded from $L^p(\rn)$ to $L^q(\rn)$ where $1<p<n/\beta$ and
$1/p-1/q=\beta/n$ (see also Paluszy\'nski \cite{p}).

For a locally integrable function $f$, the Hardy-Littlewood maximal
function $M$ is given by
$$M(f)(x)=\sup_{Q\ni x} \frac{1}{|Q|} \int_{Q} |f(y)| dy,
$$
the maximal commutator of $M$ with a locally integrable function $b$
is defined by
$$M_b(f)(x)=\sup_{Q\ni x} \frac{1}{|Q|} \int_{Q} |b(x)-b(y)||f(y)|dy,
$$
where the supremum is taken over all cubes $Q\subset\rn$ containing
$x$.

The mapping properties of the maximal commutator $M_b$ have been
studied intensively by many authors. See \cite{agkm}, \cite{ghst},
\cite{hly}, \cite{hy}, \cite{st1}, \cite{st2} and \cite{zhw2} for
instance. The following result is proved by Garc\'ia-Cuerva et al.
\cite{ghst}. See
also \cite{st1} and \cite{st2}.

\begin{thm}[\cite{ghst}] \label{thm.A} 
Let $b$ be a locally integrable function and $1<p<\infty$. Then the
maximal commutator $M_b$ is bounded from $L^p(\rn)$ to $L^p(\rn)$
if and only if $b\in{BMO(\rn)}.$
\end{thm}

The first part of this paper is to study the boundedness of $M_b$
when the symbol $b$ belongs to Lipschitz space. Some characterizations
of Lipschitz space via such commutator are given.

\begin{definition}   \label{lip}
Let $0<\beta<1$, we say a function $b$ belongs to the Lipschitz
space $\dot{\Lambda}_{\beta}(\rn)$ if there exists a constant $C$
such that for all $x,y\in\rn$,
$$|b(x)-b(y)|\le {C}|x-y|^{\beta}.
$$
The smallest such constant $C$ is called the $\dot{\Lambda}_{\beta}$
norm of $b$ and is denoted by $\|b\|_{\dot{\Lambda}_{\beta}}$.
\end{definition}

Our first result can be stated as follows.

\begin{theorem} \label{thm.mc-lp} 
Let $b$ be a locally integrable function and $0<\beta<1$,
then the following statements are equivalent:

(1)\ \ $b\in {\dot{\Lambda}_{\beta}(\rn)}$.

(2)\ \ $M_b$ is bounded from $L^p(\rn)$ to $L^q(\rn)$ for all $p,q$
with $1<p<n/\beta$ and $1/q=1/p-{\beta}/n$.

(3)\ \ $M_b$ is bounded from $L^p(\rn)$ to $L^q(\rn)$ for some $p,q$
with $1<p<n/\beta$ and $1/q=1/p-{\beta}/n$.

(4)\ \ $M_b$ satisfies the weak-type $(1,{n}/{(n-\beta)})$
estimates, namely, there exists a positive constant $C$ such that
for all $\lambda>0$,
\begin{equation}     \label{equ.weak}
\big|\{x\in\rn: M_b(f)(x)>\lambda\}\big| \le {C}
\big(\lambda^{-1}\|f\|_{L^1(\rn)}\big)^{{n}/{(n-\beta)}}.
\end{equation}

(5)\ \ $M_b$ is bounded from $L^{n/\beta}(\rn)$ to $L^{\infty}(\rn)$.
\end{theorem}

Morrey spaces were originally introduced by Morrey in \cite{m} to
study the local behavior of solutions to second order
elliptic partial differential equations. Many classical operators of
harmonic analysis were studied in Morrey type spaces during the last
decades. We refer the readers to Adams \cite{a2} and references therein.

\begin{definition}   \label{morrey}
Let $1\le {p}<\infty$ and $0\le \lambda \le {n}$. The classical
Morrey space is defined by
$$L^{p,\lambda}(\rn) = \big\{f\in {L^p_{\rm{loc}}(\rn)}:
\|f\|_{L^{p,\lambda}}<\infty \big\},
$$
where
$$\|f\|_{L^{p,\lambda}} :=\sup_{Q} \bigg(\frac1{|Q|^{\lambda/n}}
\int_Q |f(x)|^p dx\bigg)^{1/p}.
$$
\end{definition}

It is well known that if $1\le {p}<\infty$ then
$L^{p,0}(\rn)=L^p(\rn)$ and $L^{p,n}(\rn)=L^{\infty}(\rn)$.

\begin{theorem} \label{thm.mc-morrey1} 
Let $b$ be a locally integrable function and $0<\beta<1$. Suppose
that $1<p<n/\beta$, $0<\lambda<n-{\beta}p$ and
$1/q=1/p-\beta/(n-\lambda)$. Then $b \in
{\dot{\Lambda}_{\beta}}(\rn)$ if and only if $M_b$ is bounded from
$L^{p,\lambda}(\rn)$ to $L^{q,\lambda}(\rn)$.
\end{theorem}

\begin{theorem} \label{thm.mc-morrey2}
Let $b$ be a locally integrable function and $0<\beta<1$. Suppose
that $1<p<n/\beta$, $0<\lambda<n-{\beta}p$, $1/q=1/p-\beta/n$ and
$\lambda/p=\mu/q$. Then $b \in {\dot{\Lambda}_{\beta}}(\rn)$ if and
only if $M_b$ is bounded from $L^{p,\lambda}(\rn)$ to $L^{q,\mu}(\rn)$.
\end{theorem}

On the other hand, similar to (\ref{equ.1.1}), we can define
the (nonlinear) commutator of the Hardy-Littlewood maximal
function $M$ with a locally integrable function $b$ by
$$[b,M](f)(x)=b(x)M(f)(x)-M(bf)(x).
$$

Using real interpolation techniques, Milman and Schonbek \cite{ms}
established a commutator result. As an application, they obtained
the $L^p-$boundedness of $[b,M]$ when $b \in {BMO(\rn)}$ and $b\ge 0$.
This operator can be used in studying the product of a function
in $H^1$ and a function in $BMO$ (see \cite{bijz} for instance).
In 2000, Bastero, Milman and Ruiz \cite{bmr}
studied the necessary and sufficient conditions
for the boundedness of $[b,M]$ on $L^p$ spaces when $1<p<\infty$.
Zhang and Wu obtained similar results for the fractional maximal
function in \cite{zhw1} and extended the mentioned results to variable
exponent Lebesgue spaces in \cite{zhw2} and \cite{zhw3}. Recently,
Agcayazi et al. \cite{agkm} gave the end-point estimates for the
commutator $[b,M]$. Zhang \cite{zh} extended these results to the
multilinear setting.

We would like to remark that operators $M_b$ and $[b,M]$ essentially
differ from each other. For example, $M_b$ is positive and
sublinear, but $[b,M]$ is neither positive nor sublinear.

The second part of this paper is to study the mapping properties
of the (nonlinear) commutator $[b,M]$ when $b$ belongs to some
Lipschitz space.
To state our results, we recall the definition of the maximal
operator with respect to a cube. For a fixed cube $Q_0$, the
Hardy-Littlewood maximal function with respect to $Q_0$ of a
function $f$ is given by
$$M_{Q_0}(f)(x) =\sup_{Q_0\supseteq{Q}\ni{x}}
\frac{1}{|Q_0|}\int_{Q}|f(y)|dy,
$$
where the supremum is taken over all the cubes $Q$ with $Q\subseteq
{Q_0}$ and $Q\ni {x}$.

\begin{theorem}  \label{thm.nc1}
Let $b$ be a locally integrable function and $0<\beta<1$.
Suppose that $1<p<n/\beta$ and $1/q=1/p-{\beta}/n$. Then the following
statements are equivalent:

(1)\ \ $b\in {\dot{\Lambda}_{\beta}(\rn)}$ and $b\ge 0$.

(2)\ \ $[b,M]$ is bounded from $L^p(\rn)$ to $L^q(\rn)$.

(3)\ \ There exists a constant $C>0$ such that
\begin{equation}         \label{equ.nc1}
\sup_{Q} \frac1{|Q|^{\beta/n}} \bigg(\frac1{|Q|} \int_Q
|b(x)-M_{Q}(b)(x)|^q dx \bigg)^{1/q} \le {C}.
\end{equation}
\end{theorem}

\begin{theorem} \label{thm.nc2}
Let $b\ge 0$ be a locally integrable function, $0<\beta<1$ and
$b\in {\dot{\Lambda}_{\beta}(\rn)}$. Then there is a positive
constant $C$ such that for all $\lambda>0$,
$$\big|\{x\in\rn: |[b,M](f)(x)|>\lambda\}\big| \le {C}
\big(\lambda^{-1}\|f\|_{L^1(\rn)}\big)^{n/(n-\beta)}.
$$
\end{theorem}

\begin{theorem} \label{thm.nc-morrey1}
Let $b$ be a locally integrable function and $0<\beta<1$.
Suppose that $1<p<n/\beta$, $0<\lambda<n-{\beta}p$ and
$1/q=1/p-\beta/(n-\lambda)$. Then the following statements
are equivalent:

(1)\ \ $b\in {\dot{\Lambda}_{\beta}(\rn)}$ and $b\ge 0$.

(2)\ \ $[b,M]$ is bounded from $L^{p,\lambda}(\rn)$ to $L^{q,\lambda}(\rn)$.
\end{theorem}

\begin{theorem} \label{thm.nc-morrey2}
Let $b$ be a locally integrable function and $0<\beta<1$.
Suppose that $1<p<n/\beta$, $0<\lambda<n-{\beta}p$, $1/q=1/p-\beta/n$ and
$\lambda/p=\mu/q$.  Then the following statements are equivalent:

(1)\ \ $b\in {\dot{\Lambda}_{\beta}(\rn)}$ and $b\ge 0$.

(2)\ \ $[b,M]$ is bounded from $L^{p,\lambda}(\rn)$ to $L^{q,\mu}(\rn)$.
\end{theorem}

This paper is organized as follows. In the next section, we recall
some basic definitions and known results. In Section 3, we will prove
Theorems \ref{thm.mc-lp}---\ref{thm.mc-morrey2}.
Section 4 is devoted to proving Theorems \ref{thm.nc1}---\ref{thm.nc-morrey2}.

\section{Preliminaries and Lemmas}

For a measurable set $E$, we denote
by $|E|$ the Lebesgue measure and by $\chi_E$ the characteristic
function of $E$. For $p\in[1,\infty]$, we denote by $p'$ the
conjugate index of $p$, namely, $p'=p/(p-1)$. For a locally
integrable function $f$ and a cube $Q$,
we denote by $f_Q =(f)_Q =\frac{1}{|Q|} \int_{Q} f(x) dx.$

To prove the theorems, we need some known results. It is known
that the Lipschitz space $\dot{\Lambda}_{\beta}(\rn)$ coincides with
some Morrey-Companato space (see \cite{jtw} for example) and can be
characterized by mean oscillation as the following lemma, which is
due to DeVore and Sharpley \cite{ds} and Janson, Taibleson and Weiss
\cite{jtw} (see also Paluszy\'nski \cite{p}).

\begin{lemma}    \label{lem.lip}
Let $0<\beta<1$ and $1\le {q}<\infty$. Define
$$\dot{\Lambda}_{\beta,q}(\rn):=\bigg\{f \in{L_{\rm{loc}}^1(\rn)}:
\|f\|_{\dot{\Lambda}_{\beta,q}} = \sup_Q \frac1{|Q|^{\beta/n}}
\bigg(\frac1{|Q|} \int_Q|f(x)-f_Q|^qdx\bigg)^{1/q}<\infty \bigg\}.
$$
Then, for all $0<\beta<1$ and $1\le {q}<\infty$,
$\dot{\Lambda}_{\beta}(\rn)=\dot{\Lambda}_{\beta,q}(\rn)$
with equivalent norms.
\end{lemma}

Let $0<\alpha<n$ and $f$ be a locally integrable function, the
fractional maximal function of $f$ is given by
$$ \mathfrak{M}_{\alpha}(f)(x) =\sup_Q\frac1{|Q|^{1-\alpha/n}} \int_Q|f(y)|dy
$$
where the supremum is taken over all cubes $Q\subset\rn$ containing $x$.

The following strong and weak-type boundedness of $\mathfrak{M}_{\alpha}$ are
well-known, see \cite{g} and \cite{d}.

\begin{lemma}    \label{lem.frac1}
Let $0<\alpha<n$, $1\le {p}\le{n/\alpha}$ and $1/q=1/p-\alpha/n$.

(1) If $1<p<n/\alpha$ then there exists a positive constant $C(n,\alpha,p)$
such that
$$\|\mathfrak{M}_{\alpha}(f)\|_{L^q(\rn)}
\le {C(n,\alpha,p)}\|f\|_{L^p(\rn)}.
$$

(2) If $p=n/\alpha$ then there exists a positive constant  $C(n,\alpha)$
such that
$$\|\mathfrak{M}_{\alpha}(f)\|_{L^{\infty}(\rn)} \le
{C(n,\alpha)}\|f\|_{L^{n/\alpha}(\rn)}.
$$

(3) If $p=1$ then there exists a positive constant  $C(n,\alpha)$
such that for all $\lambda>0$
$$\big|\big\{x\in\rn: \mathfrak{M}_{\alpha}(f)(x)>\lambda\big\}\big|
\le {C(n,\alpha)}\big(\lambda^{-1}\|f\|_{L^1(\rn)}\big)^{n/(n-\alpha)}.
$$
\end{lemma}

Spanne (see \cite{pe}) and Adams \cite{a1} studied the boundedness
of the fractional integral $I_{\alpha}$ in classical Morrey spaces.
We note that the fractional maximal function enjoys the same
boundedness as that of the fractional integral since the pointwise
inequality $\mathfrak{M}_{\alpha}(f)(x) \le {I_{\alpha}}(|f|)(x)$.
These results can be summarized as follows (see also \cite{sh}):

\begin{lemma}  \label{lem.frac2}
Let $0<\alpha<n$, $1<p<n/\alpha$ and $0<\lambda<n-{\alpha}p$.

{\rm (1)} If $1/q=1/p-\alpha/(n-\lambda)$ then there is a
constant $C>0$ such that
$$\|\mathfrak{M}_{\alpha}(f)\|_{L^{q,\lambda}(\rn)} \le
{C}\|f\|_{L^{p,\lambda}(\rn)} ~\hbox{ for  every~} f \in
{L^{p,\lambda}(\rn)}.
$$

{\rm (2)} If $1/q=1/p-\alpha/n$ and $\lambda/p=\mu/q$. Then there
is a constant $C>0$ such that
$$\|\mathfrak{M}_{\alpha}(f)\|_{L^{q,\mu}(\rn)} \le
{C}\|f\|_{L^{p,\lambda}(\rn)} ~\hbox{ for  every~} f \in
{L^{p,\lambda}(\rn)}.
$$
\end{lemma}

\begin{lemma}[\cite{km}]    \label{lem.cube}
Let $1\le {p}<\infty$ and $0<\lambda<n$, then there is a constant $C>0$
that depends only on $n$ such that
$$\|\chi_Q\|_{L^{p,\lambda}(\rn)}\le {C}|Q|^{\frac{n-\lambda}{np}}.
$$
\end{lemma}

\section{Proof of Theorems \ref{thm.mc-lp}---\ref{thm.mc-morrey2} }

\begin{proof} {\bf of Theorem \ref{thm.mc-lp}} ~If
$b\in {\dot{\Lambda}_{\beta}(\rn)}$, then
\begin{equation}        \label{equ.3.1}
\begin{split}
M_b(f)(x) & =\sup_{Q\ni {x}} \frac1{|Q|} \int_Q
|b(x)-b(y)||f(y)|dy\\
  &\le {C}\|b\|_{\dot{\Lambda}_{\beta}} \sup_{Q\ni {x}}
    \frac1{|Q|^{1-\beta/n}} \int_Q |f(y)|dy\\
    &= {C}\|b\|_{\dot{\Lambda}_{\beta}} \mathfrak{M}_{\beta}(f)(x).
\end{split}
\end{equation}
Obviously, (2), (3), (4) and (5) follow from Lemma \ref{lem.frac1},
Lemma \ref{lem.frac2} and (\ref{equ.3.1}).

(3) $\Longrightarrow$ (1):~ Assume $M_b$ is bounded
from $L^p(\rn)$ to $L^q(\rn)$ for some $p,q$ with $1<p<n/\beta$ and
$1/q=1/p-{\beta}/n$. For any cube $Q\subset \rn$, by
H\"older's inequality and noting that $1/p+1/q'=1+\beta/n$, one gets
\begin{equation*}
\begin{split}
\frac1{|Q|^{1+\beta/n}}\int_Q|b(x)-b_Q|dx
 &\le \frac1{|Q|^{1+\beta/n}}\int_Q
    \bigg(\frac1{|Q|}\int_Q|b(x)-b(y)|dy\bigg)dx\\
 &=\frac1{|Q|^{1+\beta/n}}\int_Q
    \bigg(\frac1{|Q|}\int_Q|b(x)-b(y)|\chi_Q(y)dy\bigg)dx\\
 &\le \frac1{|Q|^{1+\beta/n}}\int_Q M_b(\chi_Q)(x)dx\\
 &\le \frac1{|Q|^{1+\beta/n}}
  \bigg(\int_Q[M_b(\chi_Q)(x)]^qdx\bigg)^{1/q}
   \bigg(\int_Q \chi_Q(x)dx\bigg)^{1/q'}\\
 &\le \frac{C}{|Q|^{1+\beta/n}} \|M_b\|_{L^p\rightarrow{L^q}}
   \|\chi_Q\|_{L^p(\rn)}\|\chi_Q\|_{L^{q'}(\rn)}\\
 &\le {C}\|M_b\|_{L^p\rightarrow{L^q}}.
\end{split}
\end{equation*}
This together with Lemma \ref{lem.lip} gives $b\in
{\dot{\Lambda}_{\beta}(\rn)}$.

(4) $\Longrightarrow$ (1):~  We assume (\ref{equ.weak})
is true and will verify
$b\in {\dot{\Lambda}_{\beta}(\rn)}$. For any fixed cube $Q_0
\subset \rn$, since for any $x\in {Q_0}$,
$$|b(x)-b_{Q_0}| \le \frac1{|Q_0|} \int_{Q_0}|b(x)-b(y)|dy,
$$
then, for all $x\in {Q_0}$,
\begin{equation*}
\begin{split}
M_b(\chi_{Q_0})(x) &=
\sup_{Q\ni{x}}\frac1{|Q|}\int_{Q}|b(x)-b(y)|\chi_{Q_0}(y)dy\\
 & \ge\frac1{|Q_0|} \int_{Q_0}|b(x)-b(y)|\chi_{Q_0}(y)dy\\
 & =\frac1{|Q_0|} \int_{Q_0}|b(x)-b(y)|dy\\
 &\ge |b(x)-b_{Q_0}|.
\end{split}
\end{equation*}

This together with (\ref{equ.weak}) gives
\begin{equation*}
\begin{split}
\big|\big\{x\in{Q_0}: |b(x)-b_{Q_0}|>\lambda \big\}\big|
&\le\big|\big\{x\in{Q_0}: M_b(\chi_{Q_0})(x)>\lambda\big\}\big|\\
 &\le{C}\big(\lambda^{-1}\|\chi_{Q_0}\|_{L^1(\rn)}\big)^{n/(n-\beta)}\\
  &={C}\big(\lambda^{-1}|Q_0|\big)^{n/(n-\beta)}.
\end{split}
\end{equation*}

Let $t>0$ be a constant to be determined later, then
\begin{equation*}
\begin{split}
\int_{Q_0}|b(x)-b_{Q_0}|dx &=\int_0^{\infty} \big|\big\{x\in{Q_0}:
|b(x)-b_{Q_0}|>\lambda \big\}\big|d\lambda\\
 & =\int_0^t \big|\big\{x\in{Q_0}:
  |b(x)-b_{Q_0}|>\lambda \big\}\big|d\lambda\\
 & \qquad + \int_t^{\infty} \big|\big\{x\in{Q_0}:
  |b(x)-b_{Q_0}|>\lambda \big\}\big|d\lambda\\
 & \le {t|Q_0|} + C\int_t^{\infty}
  \big(\lambda^{-1}|Q_0|\big)^{n/(n-\beta)} d\lambda\\
 & \le {t|Q_0|} +C|Q_0|^{n/(n-\beta)}
    \int_t^{\infty}\lambda^{-n/(n-\beta)}d\lambda\\
 &\le{C(n,\beta)}\big(t|Q_0|
 +|Q_0|^{n/(n-\beta)}t^{1-n/(n-\beta)}\big).
\end{split}
\end{equation*}
Set $t=|Q_0|^{\beta/n}$ in the above estimate, we have
$$\int_{Q_0}|b(x)-b_{Q_0}|dx \le{C}|Q_0|^{1+\beta/n}.
$$

It follows from Lemma \ref{lem.lip} that
$b\in {\dot{\Lambda}_{\beta}(\rn)}$ since $Q_0$ is
an arbitrary cube in $\rn$.

(5) $\Longrightarrow$ (1):~  If $M_b$ is bounded from $L^{n/\beta}(\rn)$ to
$L^{\infty}(\rn)$, then for any cube $Q\subset \rn$,
\begin{align*}
\frac1{|Q|^{1+\beta/n}}\int_{Q}|b(x)-b_{Q}|dx&\le\frac1{|Q|^{1+\beta/n}}
 \int_Q \bigg(\frac1{|Q|}\int_Q|b(x)-b(y)|\chi_Q(y)dy\bigg)dx\\
  & \le \frac1{|Q|^{1+\beta/n}} \int_Q M_b(\chi_Q)(x)dx\\
   & \le \frac1{|Q|^{\beta/n}}\|M_b(\chi_Q)\|_{L^{\infty}(\rn)}\\
 & \le\frac{C}{|Q|^{\beta/n}}\|M_b\|_{L^{n/\beta}\to{L^{\infty}}}
       \|\chi_Q\|_{L^{n/\beta}(\rn)}\\
 &\le {C}{\|M_b\|_{L^{n/\beta}\to{L^{\infty}}}}.
\end{align*}
This together with Lemma \ref{lem.lip} gives $b\in
{\dot{\Lambda}_{\beta}(\rn)}$.

The proof of Theorem \ref{thm.mc-lp} is completed since
(2) $\Longrightarrow$ (1) follows from (3) $\Longrightarrow$ (1).
\end{proof}

\begin{proof} {\bf of Theorem \ref{thm.mc-morrey1}} ~Assume
$b\in{\dot{\Lambda}_{\beta}(\rn)}$. By (\ref{equ.3.1}) and Lemma
\ref{lem.frac2} (1), we have
$$\|M_b(f)\|_{L^{q,\lambda}} \le\|b\|_{\dot{\Lambda}_{\beta}}
     \|\mathfrak{M}_{\beta}(f)\|_{L^{q,\lambda}}
 \le{C}\|b\|_{\dot{\Lambda}_{\beta}}\|f\|_{L^{p,\lambda}}.
$$

Conversely, if $M_b$ is bounded from $L^{p,\lambda}(\rn)$ to
$L^{q,\lambda}(\rn)$, then for any cube $Q\subset\rn$,
\begin{equation*}
\begin{split}
\frac1{|Q|^{\beta/n}}\bigg(\frac1{|Q|}\int_{Q}
 |b(x)-b_{Q}|^qdx\bigg)^{1/q}
 &\le \frac1{|Q|^{\beta/n}}\bigg(\frac1{|Q|}
 \int_Q \bigg[\frac1{|Q|}\int_Q
 |b(x)-b(y)|\chi_Q(y)dy\bigg]^q dx\bigg)^{1/q}\\
  & \le \frac1{|Q|^{\beta/n}} \bigg(\frac1{|Q|}
   \int_Q[M_b(\chi_Q)(x)]^qdx\bigg)^{1/q}\\
&=\frac1{|Q|^{\beta/n}}\bigg(\frac{|Q|^{\lambda/n}}{|Q|}\bigg)^{1/q}
 \bigg(\frac1{|Q|^{\lambda/n}}\int_Q[M_b(\chi_Q)(x)]^qdx\bigg)^{1/q}\\
&\le |Q|^{-\beta/n-1/q+\lambda/(nq)}
   \|M_b(\chi_Q)\|_{L^{q,\lambda}(\rn)}\\
&\le {C}|Q|^{-\beta/n-1/q+\lambda/(nq)}
 \|M_b\|_{L^{p,\lambda}\to{L^{q,\lambda}}}
   \|\chi_Q\|_{L^{p,\lambda}(\rn)}\\
&\le {C}\|M_b\|_{L^{p,\lambda}\to{L^{q,\lambda}}},
\end{split}
\end{equation*}
where in the last step we have used $1/q=1/p-\beta/(n-\lambda)$
and Lemma \ref{lem.cube}.

It follows from Lemma \ref{lem.lip} that $b\in{\dot{\Lambda}_{\beta}(\rn)}$.
This completes the proof.
\end{proof}

\begin{proof} {\bf of Theorem \ref{thm.mc-morrey2}} ~By a similar proof
to the one of Theorem \ref{thm.mc-morrey1}, we can obtain Theorem
\ref{thm.mc-morrey2}.
\end{proof}

\section{Proof of Theorems \ref{thm.nc1}---\ref{thm.nc-morrey2}}

\begin{proof}{\bf of Theorem \ref{thm.nc1}} ~(1) $\Longrightarrow$ (2):~
For any fixed $x\in\rn$ such that $M(f)(x)<\infty$, since $b\ge0$ then
\begin{equation}   \label{equ.proof-nc1-1}
\begin{split}
|[b,M](f)(x)|&=|b(x)M(f)(x)-M(bf)(x)|\\
 & =\bigg|\sup_{Q\ni{x}} \frac1{|Q|} \int_Q b(x)|f(y)|dy
    -\sup_{Q\ni{x}} \frac1{|Q|} \int_Q b(y)|f(y)|dy\bigg|\\
 & \le \sup_{Q\ni{x}} \frac1{|Q|} \int_Q |b(x)-b(y)||f(y)|dy\\
 & = M_b(f)(x).
\end{split}
\end{equation}
It follows from Theorem \ref{thm.mc-lp} that $[b,M]$ is
bounded from $L^p(\rn)$ to $L^q(\rn)$ since
$b\in{\dot{\Lambda}_{\beta}(\rn)}$.

(2) $\Longrightarrow$ (3):~ For any fixed cube $Q\subset\rn$
and all $x\in {Q}$, we have (see the proof of Proposition 4.1 in
\cite{bmr}, see also (2.4) in \cite{zhw1})
$$M(\chi_Q)(x)=\chi_Q(x) ~~\hbox{and}~~ M(b\chi_Q)(x)=M_Q(b)(x).
$$

Then,
\begin{equation}              \label{equ.proof-nc1-2}
\begin{split}
&\frac1{|Q|^{\beta/n}} \bigg(\frac1{|Q|}\int_Q
   \big|b(x)-M_Q(b)(x)\big|^qdx\bigg)^{1/q}\\
 &=\frac1{|Q|^{\beta/n}} \bigg(\frac1{|Q|}\int_Q
   \big|b(x)M(\chi_Q)(x)-M_Q(b\chi_Q)(x)\big|^qdx\bigg)^{1/q}\\
  &=\frac1{|Q|^{\beta/n}} \bigg(\frac1{|Q|}\int_Q
    \big|[b,M](\chi_Q)(x)\big|^qdx\bigg)^{1/q}\\
  &\le \frac1{|Q|^{\beta/n+1/q}} \big\|[b,M](\chi_Q)\big\|_{L^q(\rn)}\\
    &\le \frac{C}{|Q|^{\beta/n+1/q}} \big\|\chi_Q\big\|_{L^p(\rn)}\\
    &\le{C},
\end{split}
\end{equation}
which implies (3) since the cube $Q\subset \rn$ is arbitrary.

(3) $\Longrightarrow$ (1):~ To prove $b\in{\dot{\Lambda}_{\beta}(\rn)}$,
by Lemma \ref{lem.lip}, it suffices to verify that there is a
constant $C>0$ such that for all cubes $Q$,
\begin{equation}              \label{equ.proof-nc1-3}
\begin{split}
\frac1{|Q|^{1+\beta/n}} \int_Q |b(x)-b_Q|dx \le {C}.
\end{split}
\end{equation}

For any fixed cube $Q$, let $E=\{x\in {Q}: b(x)\le {b_Q}\}$ and
$F=\{x\in {Q}: b(x)>b_Q\}$. The following equality is
trivially true (see \cite{bmr} page 3331):
$$\int_E|b(x)-b_Q|dx=\int_F|b(x)-b_Q|dx.
$$

Since for any $x\in{E}$ we have $b(x)\le {b_Q} \le {M_Q(b)(x)}$, then
for any  $x\in{E}$,
$$|b(x)-b_Q| \le |b(x)-M_Q(b)(x)|.
$$
Thus,
\begin{equation}              \label{equ.proof-nc1-4}
\begin{split}
\frac1{|Q|^{1+\beta/n}} \int_Q |b(x)-b_Q|dx
  &=\frac1{|Q|^{1+\beta/n}}\int_{E\cup{F}}|b(x)-b_Q|dx\\
  &=\frac2{|Q|^{1+\beta/n}}\int_E|b(x)-b_Q|dx\\
    &\le \frac2{|Q|^{1+\beta/n}}\int_E|b(x)-M_Q(b)(x)|dx\\
        &\le \frac2{|Q|^{1+\beta/n}}\int_Q|b(x)-M_Q(b)(x)|dx.
\end{split}
\end{equation}

On the other hand, it follows from
H\"older's inequality and (\ref{equ.nc1}) that
\begin{equation*}
\begin{split}
\frac1{|Q|^{1+\beta/n}}\int_Q\big|b(x)-M_Q(b)(x)\big|dx &
 \le \frac1{|Q|^{1+\beta/n}}\bigg(\int_Q\big|b(x)-M_Q(b)(x) \big|^qdx\bigg)^{1/q}
  |Q|^{1/q'} \\
  &\le \frac1{|Q|^{\beta/n}}\bigg(\frac1{|Q|}
    \int_Q\big|b(x)-M_Q(b)(x) \big|^qdx\bigg)^{1/q}\\
    &\le {C}.
\end{split}
\end{equation*}
this together with (\ref{equ.proof-nc1-4}) gives (\ref{equ.proof-nc1-3})
and so we achieve $b\in{\dot{\Lambda}_{\beta}(\rn)}$.

In order to prove $b\ge 0$, it suffices to show $b^{-}=0$, where
$b^{-}=-\min\{b,0\}$. Let $b^{+}=|b|-b^{-}$, then
$b=b^{+}-b^{-}$. For any fixed cube $Q$, observe that
$$0\le {b^{+}(x)}\le |b(x)| \le {M_Q(b)(x)}$$
for $x\in{Q}$ and therefore we have that, for $x\in{Q}$,
$$0\le {b^{-}(x)}\le {M_Q(b)(x)}- b^{+}(x) + {b^{-}(x)}
= M_Q(b)(x)- b(x).
$$

Then, it follows from (\ref{equ.nc1}) that, for any cube $Q$,
\begin{align*}
& \frac1{|Q|}\int_Q b^{-}(x)dx \le \frac1{|Q|}\int_Q |M_Q(b)(x)- b(x)|\\
&\le \bigg(\frac1{|Q|}\int_Q |b(x)-M_Q(b)(x)|^q dx\bigg)^{1/q}\\
&= |Q|^{\beta/n} \bigg\{ \frac1{|Q|^{\beta/n}}
 \bigg(\frac1{|Q|}\int_Q |b(x)-M_Q(b)(x)|^q dx\bigg)^{1/q}\bigg\}\\
&\le {C}|Q|^{\beta/n}.
\end{align*}
Thus, $b^{-}=0$ follows from Lebesgue's differentiation theorem.

The proof of Theorem \ref{thm.nc1} is completed.
\end{proof}

\begin{proof}{\bf of Theorem \ref{thm.nc2}} ~Obviously, Theorem \ref{thm.nc2}
follows from (\ref{equ.proof-nc1-1}) and Theorem \ref{thm.mc-lp}.
\end{proof}

\begin{proof}{\bf of Theorem \ref{thm.nc-morrey1}} ~(1) $\Longrightarrow$ (2):~
Assume $b\ge 0$ and
$b\in{\dot{\Lambda}_{\beta}(\rn)}$, then by (\ref{equ.proof-nc1-1}) and
Theorem \ref{thm.mc-morrey1} we see that $[b,M]$ is bounded from
$L^{p,\lambda}(\rn)$ to $L^{q,\lambda}(\rn)$.

(2) $\Longrightarrow$ (1):~  Assume that $[b,M]$ is bounded from $L^{p,\lambda}(\rn)$ to
$L^{q,\lambda}(\rn)$. Similarly to (\ref{equ.proof-nc1-2}), we have, for
any cube $Q\subset \rn$,
\begin{align*}
&\frac1{|Q|^{\beta/n}} \bigg(\frac1{|Q|}\int_Q
   \big|b(x)-M_Q(b)(x)\big|^qdx\bigg)^{1/q}\\
  &=\frac1{|Q|^{\beta/n}} \bigg(\frac1{|Q|}\int_Q
    \big|[b,M](\chi_Q)(x)\big|^qdx\bigg)^{1/q}\\
  &\le \frac{|Q|^{{\lambda}/(nq)}}{|Q|^{\beta/n+1/q}}
     \big\|[b,M](\chi_Q)\big\|_{L^{q,\lambda}(\rn)}\\
   &\le \frac{C|Q|^{{\lambda}/(nq)}}{|Q|^{\beta/n+1/q}}
       \big\|\chi_Q\big\|_{L^{p,\lambda}(\rn)}\\
    &\le{C},
\end{align*}
where in the last step we have used $1/q=1/p-\beta/(n-\lambda)$
and Lemma \ref{lem.cube}.

This shows by Theorem \ref{thm.nc1} that
$b\in{\dot{\Lambda}_{\beta}(\rn)}$ and $b\ge 0$.
\end{proof}

\begin{proof}{\bf of Theorem \ref{thm.nc-morrey2}} ~By the same way of
the proof of Theorem \ref{thm.nc-morrey1}, Theorem \ref{thm.nc-morrey2}
can be proven. We omit the details.
\end{proof}



\end{document}